\documentclass[a4paper]{article}

\usepackage[T1]{fontenc}
\usepackage[utf8]{inputenc}
\usepackage[english]{babel}
\usepackage{graphics}
\usepackage{epsfig}
\usepackage{amsmath}
\usepackage{amssymb}
\usepackage{amsthm}
\usepackage{url}
\usepackage{hyperref}
\usepackage{float}
\usepackage{pifont}
\usepackage[table]{xcolor}
\usepackage{soul}
\usepackage{lmodern}

\usepackage{draftwatermark} 
\SetWatermarkText{Preprint} 
\SetWatermarkScale{1}       

\DeclareMathOperator{\R}{\mathbb{R}}
\DeclareMathOperator{\re}{Re}
\DeclareMathOperator{\e}{e}
\DeclareMathOperator{\D}{\mathrm{d}\hskip-0.4ex}

\title{A Note on the Cross Gramian for Non-Symmetric Systems}
\author{Christian Himpe\thanks{Contact: \href{mailto:christian.himpe@uni-muenster.de}{\nolinkurl{christian.himpe@uni-muenster.de}}, 
                                        \href{mailto:mario.ohlberger@uni-muenster.de}{\nolinkurl{mario.ohlberger@uni-muenster.de}}, 
Institute for Computational and Applied Mathematics at the University of M\"unster, Einsteinstrasse~62, D-48149 M\"unster, Germany} \and Mario Ohlberger\footnotemark[1]}

\date{}

\newtheoremstyle{thm}{\topsep}{\topsep}{\normalfont \itshape}{}{\normalfont \bfseries}{}{\newline}{}
\theoremstyle{thm}

\newcounter{dummy}

\newtheorem{mydefine}[dummy]{Definition}

\begin{document}

\maketitle

\begin{abstract}\bfseries
The cross gramian matrix is a tool for model reduction and system identification, but it is only computable for square control systems.
For symmetric systems the cross gramian possesses a useful relation to the system's associated Hankel singular values.
Yet, many real-life models are neither square nor symmetric.
In this work, concepts from decentralized control are used to approximate a cross gramian for non-symmetric and non-square systems.
To illustrate this new non-symmetric cross gramian, it is applied in the context of model order reduction.
\end{abstract}
~

\section{Introduction}
The cross gramian was introduced in \cite{fernando83} for single-input-single-output (SISO) systems and extended to multiple-input-multiple-output (MIMO) systems in \cite{laub83}.
With many applications in model order reduction (MOR) and system theory in general, such as system identification \cite{mironovskii15}, decentralized control \cite{moaveni06}, parameter identification \cite{himpe14a} or sensitivity analysis \cite{streif06}, 
a major hindrance in the use of the cross gramian matrix is the constraint that it can be computed strictly for square systems and exhibits its core property only for symmetric systems.
An extension of the cross gramian to non-symmetric systems enables such uses and particularly the (approximate) balancing state-space reduction by the computation of a single gramian matrix.
This work can be seen as a follow up to \cite{laub83} and \cite{abreu86}, which previously expanded the scope of cross gramian.

The object of interest in this context is a linear time-invariant state-space system:
\begin{align} \label{ss}
 \dot{x}(t) &= Ax(t) + Bu(t), \notag \\
       y(t) &= Cx(t) + Du(t), \\
       x(0) &= x_0, \notag
\end{align}
which consists of a linear vector field composed of a system matrix $A \in \R^{N \times N}$ and an input matrix $B \in \R^{N \times M}$, as well as a linear output functional composed of an output matrix $C \in \R^{O \times N}$ and a feed-forward matrix $D \in \R^{O \times M}$.
In the scope of this work, only a trivial feed-forward matrix $D = 0$ is considered.
Furthermore, the system is assumed to be asymptotically stable, hence all eigenvalues of the system matrix $A$ lie in the left half-plane $\re(\lambda_i(A))<0$.
The classic cross gramian is defined for a subset of these systems and is expanded to arbitrary systems in the remainder of this work, which is organized as follows.
In \hyperref[cg]{Section~\ref*{cg}} the cross gramian is reviewed.
Existing approaches and the new result for the non-symmetric cross gramian are presented in \hyperref[ncg]{Section~\ref*{ncg}}.
Lastly, in \hyperref[num]{Section~\ref*{num}}, verification and validation for the new non-symmetric cross gramian and comparison with established methods is conducted.

\section{The Cross Gramian} \label{cg}
With the controllability operator $\mathcal{C}(u) := \int_0^\infty \e^{At}Bu(t) \D t$ and the observability operator $\mathcal{O}(x) := C\e^{At}x$, the controllability and observability of a system can be evaluated through the associated controllability gramian $W_C := \mathcal{C}\mathcal{C}^*$ and observability gramian $W_O := \mathcal{O}^*\mathcal{O}$. 
A third system gramian, called cross gramian, combines controllability and observability into a single matrix.

\subsection{Square Systems}
For a square system, a system with the same number of inputs and outputs, the cross gramian is defined as the product of controllability operator $\mathcal{C}$ and observability operator $\mathcal{O}$ \cite{fernando83}:
\begin{align} \label{wx}
 W_X := \mathcal{C}\mathcal{O} = \int_0^\infty \e^{At} BC \e^{At} \D t \in \R^{N \times N}.
\end{align}
Classically, the cross gramian is computed as solution to the Sylvester equation:
\begin{align*}
 A W_X + W_X A = -BC,
\end{align*}
which relates to the definition in \eqref{wx} through integration by parts:
\begin{align*}
 \int_0^\infty \e^{At} BC \e^{At} \D t &= A^{-1} \begin{vmatrix} \e^{At} BC \e^{At} \end{vmatrix}_0^\infty \\ & - A^{-1} \int_0^\infty \e^{At} BC \e^{At} \D t A \\
 \Rightarrow A \int_0^\infty \e^{At} BC \e^{At} \D t &=  -BC \\ &- \int_0^\infty \e^{At} BC \e^{At} dt A \\
 \Rightarrow A W_X &= -BC - W_X A.
\end{align*}

\subsection{Symmetric Systems}
A system in the form of \eqref{ss} is called symmetric if the system's gain is symmetric\footnote{Equivalently the symmetry of the impulse response, transfer function or Markov parameter can be used.}:
\begin{align*}
 CA^{-1}B = (CA^{-1}B)^T.
\end{align*}
In other words for a symmetric system exists a symmetric matrix $P$ such that:
\begin{align} \label{eq:sym}
 AP = PA^T \text{ and } B = P^{-1} C^T
\end{align}
are fulfilled \cite{willems76,fortuna88}.
Since a SISO system has a scalar gain, all SISO systems are symmetric.

For a symmetric system, the absolute value of the eigenvalues of the cross gramian are equal to the Hankel singular values \cite{sorensen02}:
\begin{align} \label{ll}
 W_X^2 &= W_C W_O \notag \\
 \Rightarrow |\lambda(W_X)| &= \sqrt{\lambda(W_C W_O)}.
\end{align}
This core property of the cross gramian allows to evaluate controllability and observability information of a system by computing a single gramian matrix.
Numerically, this can be more efficient than computing two, namely the controllability and observability gramians, which additionally have to be balanced to be used, in example, for model order reduction through balanced truncation \cite{baur14}.

\section{The Non-Symmetric Cross Gramian} \label{ncg}
In this section, existing approaches for cross gramians of non-symmetric systems and selected methods from decentralized control are briefly summarized; the latter is employed to expand the scope of the cross gramian from symmetric to non-symmetric systems.

\subsection{Previous Work} \label{alt}
To the authors' best knowledge there exist two methods to broaden the scope of the cross gramian for non-symmetric MIMO systems and towards more general configurations.

The first approach from \cite{abreu86} extends the applicability of the cross gramian from symmetric systems to the wider class of orthogonally symmetric systems.
Given a symmetric matrix $P = P^T$ for which $AP = PA^T$ holds and an orthogonal matrix $U$, with the property $B = PCU^T$ if $O \le M$, or $C = PBU^T$ if $M \le O$, then the system is orthogonally symmetric and the associated cross gramian:
\begin{align*}
 \breve{W}_X &= \int_0^\infty \e^{At}BUC \e^{At} \D t,
\end{align*}
satisfies the core property \eqref{ll}.

The second approach, presented in \cite{sorensen01,sorensen02}, uses embedding of a non-symmetric or non-square system into a symmetric system,
and relies on a symmetrizer matrix \cite{datta88} $J = J^T$:
\begin{align*}
 AJ = JA^T.
\end{align*}
For a symmetrizer matrix $J$ to $A$, an embedding is given by:
\begin{align*}
 \dot{\hat{x}}(t) &= Ax(t) + \begin{pmatrix} J C^T & B\end{pmatrix}u(t), \\
      \hat{y}(t)  &= \begin{pmatrix}C \\ B^T J^{-1}\end{pmatrix}x(t);
\end{align*}
the associated cross gramian has the form:
\begin{align}\label{aug}
 \widehat{W}_X &= \int_0^\infty \e^{At} (J C^T C + B B^T J^{-1}) \e^{At} \D t.
\end{align}
While the first approach preserves the central property from \eqref{ll} in $\breve{W}_X$ for orthogonally symmetric systems, 
the second approach approximates it in $\widehat{W}_X$ for arbitrary systems by an embedding.

\subsection{System Gramian Decomposition}
The basis for the non-symmetric cross gramian is a result from decentralized control, which aims to partition MIMO systems into sets of SISO systems with input-output pairings that exhibit the strongest coherence.
To this end a relation between the MIMO system gramians and the associated SISO subsystem gramians is described.
As a first step a MIMO system is decomposed into $M \times O$ SISO systems, by partitioning the input matrix $B$ and output matrix $C$ column-wise and row-wise respectively:
\begin{align*}
 B = \begin{pmatrix} b_1 & \hdots & b_M \end{pmatrix}, b_i \in \mathbb{R}^{N \times 1}, \; C = \begin{pmatrix} c_1 \\ \vdots \\ c_O \end{pmatrix}, c_j \in \mathbb{R}^{1 \times N}.
\end{align*}
Each combination of $b_i$ and $c_j$ induces a SISO system with the following system gramians:
\begin{align*}
 W_C^i &:= \int_0^\infty \e^{At} b_i b_i^T \e^{A^Tt} \D t, \\
 W_O^j &:= \int_0^\infty \e^{A^Tt} c^T_j c_j \e^{At} \D t, \\
 W_X^{i,j} &:= \int_0^\infty \e^{At} b_i c_j \e^{At} \D t.
\end{align*}

The system gramians computed for the SISO subsystems relate to the full MIMO gramians as shown in \cite{moaveni06}: 
\begin{align}\label{dcwcwo}
 W_C &= \sum_{i=1}^M W_C^i, \\
 W_O &= \sum_{j=1}^O W_O^j. \notag
\end{align}
For square systems, also the cross gramian can be computed and the following identity holds \cite{moaveni06}:
\begin{align}\label{dcwx}
 W_X = \sum_{i=1}^{M(=O)} W_X^{i,i}.
\end{align}

\subsection{Main Result}
Next, the previous result is utilized for the computation of an approximate cross gramian for non-square or non-symmetric systems.
The central idea is to exploit, that for any SISO system a cross gramian with the property \eqref{ll} can be computed.

With \eqref{dcwcwo}, the product of controllability and observability can be expressed as sum of squared cross gramians:
\begin{align}
 W_C W_O = \sum_{i=1}^M \sum_{j=1}^O W_C^i W_O^j = \sum_{i=1}^M \sum_{j=1}^O W_X^{i,j} W_X^{i,j}.
\end{align}
Due to the squaring of $W_X^{i,j}$, this ansatz is not numerically efficient.
Therefore, an alternative approximate non-symmetric cross gramian, related to \eqref{dcwx}, is introduced.
\begin{mydefine}[Non-Symmetric Cross Gramian]\label{nonsym}
The non-symmetric cross gramian is defined as the sum of the cross gramians of all $M \times O$ SISO subsystems:
 \begin{align} \label{cp}
  W_Z := \sum_{i=1}^M \sum_{j=1}^O W_X^{i,j}.
 \end{align}
\end{mydefine}
Obviously, this gramian does not preserve the cross gramian's property \eqref{ll}
and it should be emphasized that the non-symmetric cross gramian does not reduce to the classic cross gramian in case of a symmetric system (compare \eqref{dcwx}).
Yet, for a linear system, $W_Z$ yields the following representation:
\begin{align} \label{eq:mn}
 W_Z &= \sum_{i=1}^M \sum_{j=1}^O \int_0^\infty \e^{At} b_i c_j \e^{At} \D t \notag \\
     &= \int_0^\infty \e^{At} \sum_{i=1}^M \sum_{j=1}^O b_i c_j \e^{At} \D t \\
     &= \int_0^\infty \e^{At} (\sum_{i=1}^M b_i) (\sum_{j=1}^O c_j) \e^{At} \D t. \notag
\end{align}
Hence, this approximate cross gramian $W_Z$ is equal to the cross gramian of the SISO system given by $A$, the row sum of the input matrix $B$ and the column sum of the output matrix $C$,
and thus can be seen as an ``average'' cross gramian over all SISO subsystems.

Both approaches in \hyperref[alt]{Section~\ref*{alt}} share the common drawback, that they are limited to linear systems \eqref{ss}, and additionally may require a, potentially computationally expensive, symmetrizer.
The new approach, proposed in \eqref{cp}, does not require the linear structure of the underlying system or derive system properties using linear algebra.
Only a cross gramian of a SISO system needs to be computable.
Thus, the non-symmetric cross gramian can even be computed for nonlinear systems if a nonlinear cross gramian is available;
this could be for example an empirical cross gramian \cite{himpe14a}.
Empirical gramians \cite{lall99} are computed from (simulated) trajectories of the underlying system with perturbations in input and initial state.
The empirical cross gramian has been introduced for SISO systems in \cite{streif06} and extended to MIMO systems in \cite{himpe14a};
and as shown in \cite{himpe14a}, the empirical cross gramian is equal to the cross gramian \eqref{wx} for linear systems \eqref{ss}, hence the empirical cross gramian can be used for the following experiments.

\section{Numerical Results} \label{num}
The presented method is implemented as part of the empirical gramian framework (\texttt{emgr}) \cite{emgr} introduced in \cite{himpe13a};
and the following numerical experiments are conducted using \texttt{emgr}, which is compatible with \textsc{OCTAVE} and \textsc{MATLAB}\textsuperscript{\textregistered}.
The source code for reproducing the experiments can be found under an open source license in the supplemental materials and at \url{http://runmycode.org/companion/view/913} .

From a computational point of view, the proposed method has the advantage, that for all SISO subsystem cross gramians the trajectories for perturbed initial states (observability), which consume the dominant fraction of overall computational time, have to be computed only once.

Next, the presented non-symmetric cross gramian is tested in the context of projection-based model reduction \cite{baur14}:
For a linear system \eqref{ss} a reduced order model is obtained through projection-based model reduction using the (truncated) projection $U_1,V_1$ by:
\begin{align*}
 \dot{x}_r &= V_1 A U_1 x_r(t) + V_1 B u(t), \\
       y_r &= C U_1 x_r(t).
\end{align*}
Among others, such projections can be computed by balancing approaches, like balanced truncation \cite{moore81} utilizing the controllability gramian and observability gramian, or direct truncation (approximate balancing) utilizing the cross gramian.
By truncating the left singular vectors from the singular value decomposition of the cross gramian a (one-sided) Galerkin projection is generated, for further details see \cite{himpe14a},
\begin{align*}
 W_X \stackrel{SVD}{=} U D V \rightarrow U = \begin{pmatrix} U_1 & U_2 \end{pmatrix} \rightarrow V_1 := U_1^T.
\end{align*}
Generally, a Galerkin projection cannot be expected to be stability-preserving, yet in case of the cross gramian the same argument as for balanced POD in \cite[Ch.~5.4]{holmes12} can be made due to the use of an energy-preserving inner-product,
\begin{align*}
 W_X W_X^T = \mathcal{CO}(\mathcal{CO})^* = \langle \mathcal{C}, \mathcal{C}^* \rangle_{\mathcal{O}\mathcal{O}^*}.
\end{align*}
This property transfers also to the newly introduced non-symmetric cross gramian, since it reduces to a classic cross gramian for the ``averaged'' system \eqref{eq:mn} that has the same system matrix $A$.

First, the approximate non-symmetric cross gramian \eqref{cp} is compared to balanced truncation \cite{moore81,baur14} and the classic cross gramian \eqref{wx} for a symmetric system.

A state-space symmetric system $A = A^T$, $B = C^T$ of state-space dimension $N = \dim(x(t)) = 64$ and input dimension $M = \dim(u(t)) = \dim(y(t)) = O = 8$ is selected, 
with a negative Lehmer matrix as system matrix $A$, $A_{ij} := -\frac{\min(i,j)}{\max(i,j)}$ and uniformly random generated input matrix $B = C^T$.

To confirm \eqref{eq:mn}, the error between the proposed non-symmetric cross gramian $W_Z$ and the cross gramian $\overline{W}_X$ of the SISO system $\{A,b,c\}$, with $b_i = \sum_{j=1}^M B_{ij}$ and $c_j = \sum_{i=1}^O C_{ij}$, is compared in the Frobenius norm:
\begin{align*}
 \| W_Z - \overline{W}_X \|_F \approx 1.02\text{\sc{e}-}13,
\end{align*}
which is close to machine precision in double precision floating point arithmetic.

For the following experiments a zero initial state $x_0 = 0$ is set and an impulse input $u_i(t) = \delta(t), \: i=1\dots M$, is applied.
The relative $L_2$ output error is then evaluated for varying reduced order state-space dimensions.
Such an coherence assessment for a reduced order impulse response is meaningful since a linear system response is a convolution of the impulse response with an input function.

\begin{figure}
\includegraphics[width=\textwidth]{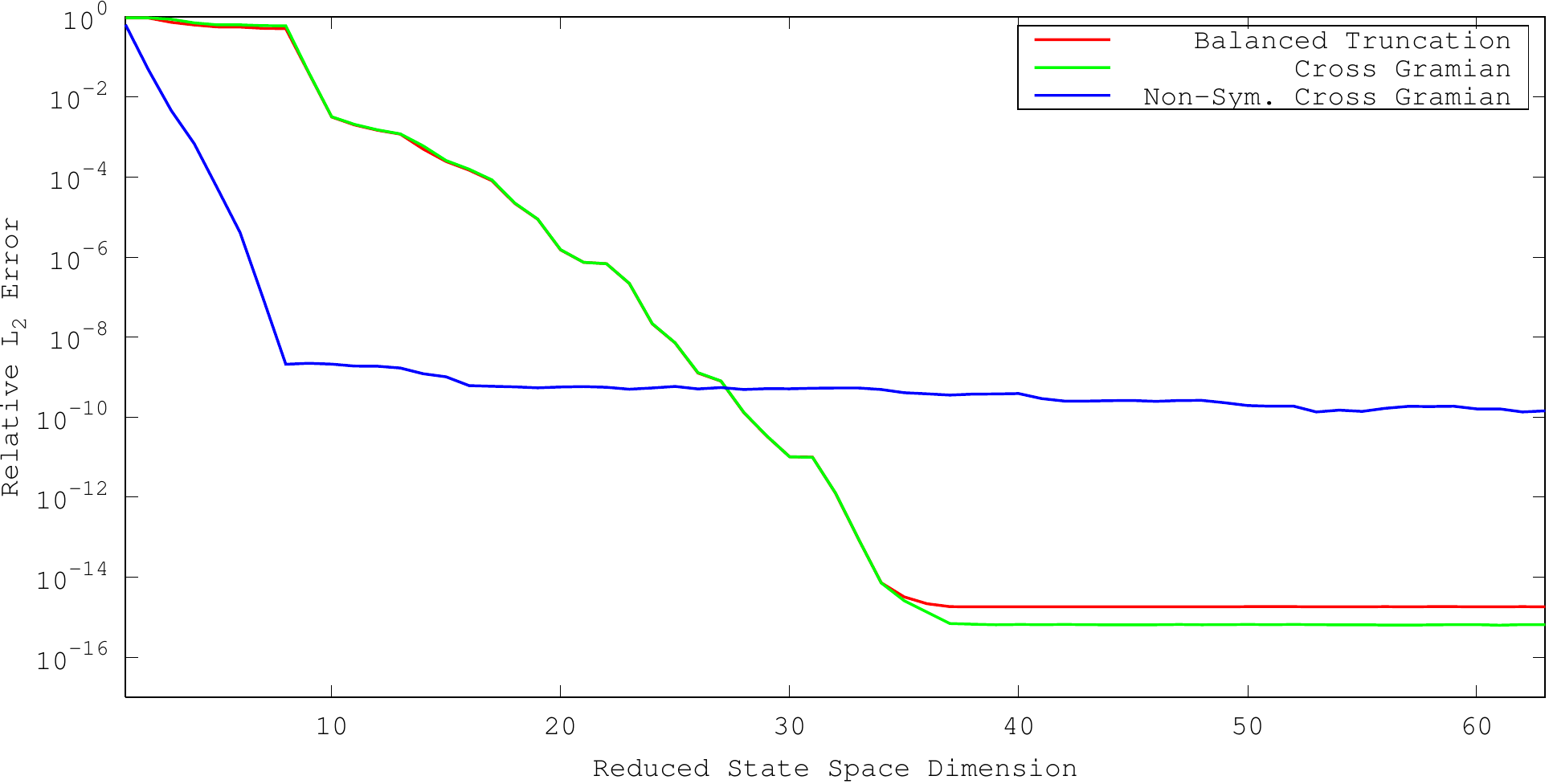}
\caption{Relative $L_2$ output error of reduced order models from balanced truncation, cross gramian and  non-symmetric cross gramian for a state-space symmetric system.}
\label{bt}
\end{figure}

\hyperref[bt]{Figure~\ref*{bt}} shows, that the reduced order models obtained by balanced truncation and the cross gramian exhibits the same behavior due to state-space symmetric nature of the system.
The newly proposed non-symmetric cross gramian from \hyperref[nonsym]{Definition~\ref*{nonsym}} does not achieve the same accuracy, but provides a lower output error for lesser order reduced models.
A lower accuracy is to be expected due to the use of a one-sided projection; while balanced truncation uses a two-sided Petrov-Galerkin projection, which in the state-space symmetric case is equal to the projection obtained from the (symmetric) classic cross gramian.  
Notably, the output error of the reduced order model from the non-symmetric cross gramian drops already at a very low order $n \geq 8$ to the steady error level while balanced truncation and the classic cross gramian reach this error not until a reduced order of $n \geq 38$. 

Second, the non-symmetric cross gramian is compared to balanced truncation and the approximate cross gramian \eqref{aug} of the symmetric system derived by embedding from \eqref{aug} for a non-square (and thus non-symmetric) system.

\begin{figure}
\includegraphics[width=\textwidth]{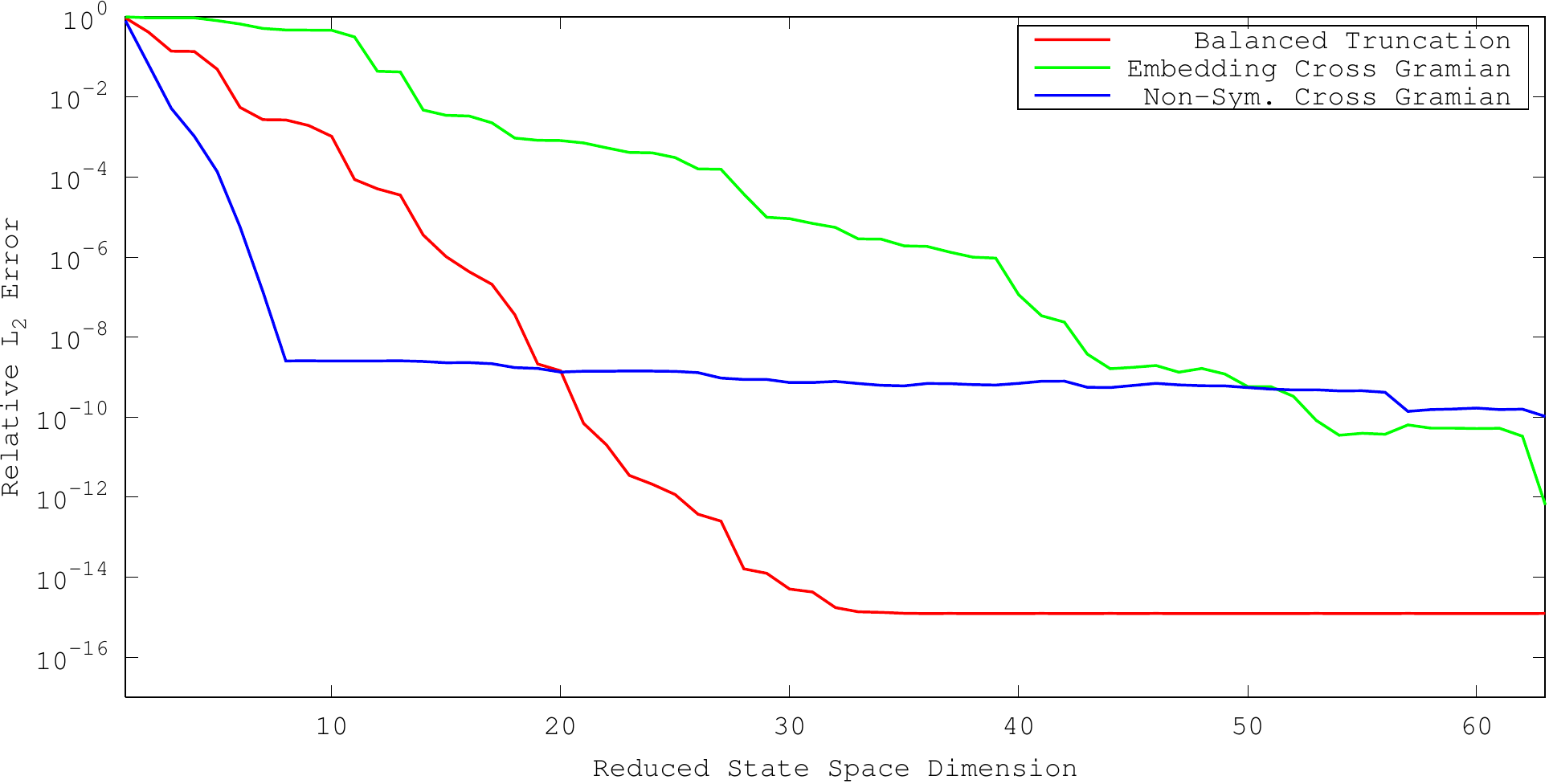}
\caption{Relative $L_2$ output error of reduced order models from balanced truncation, embedding cross gramian and non-symmetric cross gramian for a non-square system.}
\label{ct}
\end{figure}

To prevent the computation of a symmetrizer matrix $J$, the symmetric system matrix $A \in \R^{64 \times 64}$ from the first example is reused,
yet now, a uniformly random generated input matrix $B \in \R^{64 \times 4}$ and a uniformly random output matrix $C \in \R^{8 \times 64}$ is selected, thus the system is non-square and non-symmetric, since $M = 4$ and $O = 8$.
Also for this example, zero initial state $x_0 = 0$ and impulse input $u_i(t) = \delta(t),\: i=1\dots4$, is applied and the relative $L_2$ output error is evaluated for varying reduced order state-space dimensions in \hyperref[ct]{Figure~\ref*{ct}}.

With the reference of balanced truncation, the cross gramian of the embedded system performs worse as predicted in \cite{sorensen01,sorensen02}.
Again, the non-symmetric cross gramian yields reduced models with less relative output error for small (reduced) orders.

Lastly, a square but non-symmetric stable system with system matrix $A \in \R^{64 \times 64}$, uniformly random input matrix $B \in \R^{64 \times 8}$ and uniformly random output matrix $C \in \R^{8 \times 64}$ is chosen.
For this system, balanced truncation, the cross gramian and the non-symmetric cross gramian are compared.
Due to the non-symmetric nature of the system the cross gramian has no theoretical foundation, yet is applicable since the system is square.
Heuristically, also in \cite{baur08} workable results for non-symmetric systems have been obtained for the classic cross gramian. 

\begin{figure}
\includegraphics[width=\textwidth]{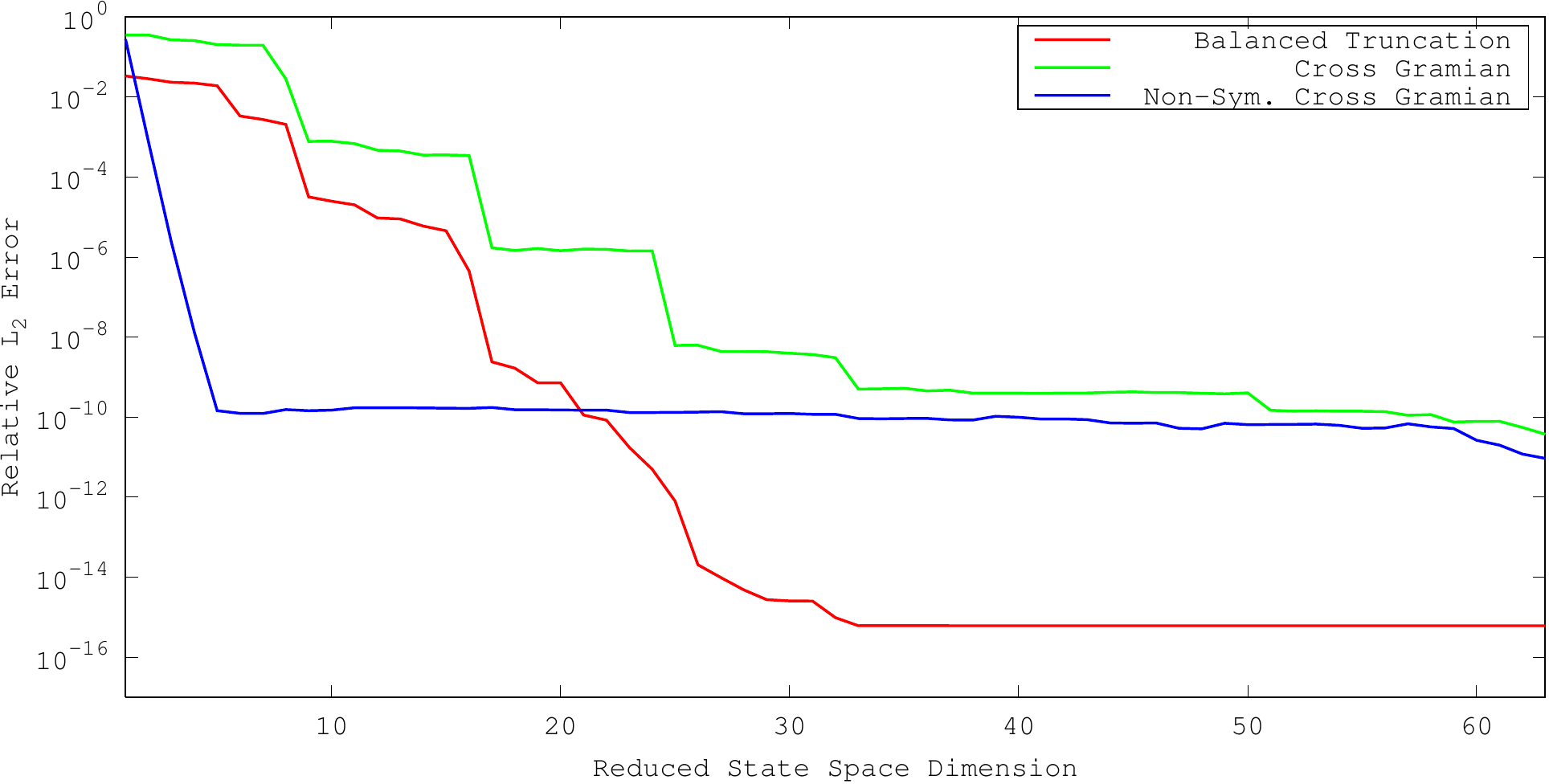}
\caption{Relative $L_2$ output error of reduced order models from balanced truncation, cross gramian and non-symmetric cross gramian for a non-symmetric system.}
\label{nt}
\end{figure}

\hyperref[nt]{Figure~\ref*{nt}} shows similar results as for the non-square system: the cross gramian and non-symmetric cross gramian both a lower accuracy,
yet the non-symmetric cross gramian requires the least states for this error level even compared to balanced truncation.

For all three experiments, the state-space model reduction error of the presented non-symmetric cross gramian is worse than for balanced truncation, but the descent of the error is steeper.
Hence, if the lower accuracy is acceptable, a smaller reduced order model can be constructed with this (non-symmetric) variant of the cross gramian method.
This superior performance of the non-symmetric cross gramian compared to balanced truncation for low order reduced order models requires further investigation.

It should be noted, that the non-symmetric cross gramian is not suitable for frequency-space-based approximation.
This is due to its construction by averaging the sub-system cross gramians.
For example, Hardy $\mathcal{H}_2$- and $\mathcal{H}_\infty$-errors between full and reduced order model will decay very slowly or not at all,
hence this method is targeted purely at state-space approximations.

\section{Conclusion}
In this work a non-symmetric cross gramian, based on concepts from decentralized control, is proposed, 
and demonstrated to provide viable results for linear non-symmetric and non-square systems, which are outside the scope of the regular cross gramian.

Future work will evaluate the effectiveness of the non-symmetric cross gramian for parametric and nonlinear systems. 
Furthermore, since the procedure suggested in \cite{rahrovani14} does not yield stable reduced order models in this setting, an investigation of alternative two-sided projections \cite{gugercin03} for the (non-symmetric) cross gramian, may yield errors comparable to balanced truncation.

\section*{Acknowledgement}
This work was supported by the Deutsche Forschungsgemeinschaft, the Open Access Publication Fund of the University of M\"unster, DFG EXC 1003 Cells in Motion - Cluster of Excellence, M\"unster, Germany 
as well as by the Center for Developing Mathematics in Interaction, DEMAIN, M\"unster, Germany.

\bibliographystyle{plainurl}
\bibliography{nonsym}







\end{document}